\documentclass{birkjour}
 \newtheorem{thm}{Theorem}[section]
 \newtheorem{cor}[thm]{Corollary}

 \newtheorem{fac}[thm]{Fact}
 \theoremstyle{definition}
 
 \theoremstyle{remark}
 \newtheorem{rem}[thm]{Remark}
 \newtheorem{rems}[thm]{Remarks}
 \newtheorem{ttt}[thm]{}
 \newtheorem{cono}[thm]{Construction}
 \newtheorem{conso}[thm]{Constructions}
 \newtheorem*{ex}{Example}
 \numberwithin{equation}{section}

\usepackage{amssymb}
\usepackage{amsmath}
\usepackage{bbm}
\usepackage{mathrsfs}

\newcounter{abc}
\newenvironment{abc}{\begin{list}{\rm \alph{abc}) }{\usecounter{abc} \leftmargin=0.0pt \labelsep=0.0pt \listparindent=0.0pt \labelwidth=0.0pt \parsep=\smallskipamount \itemsep=0.0pt \topsep=0.0pt \partopsep=\smallskipamount}}{\end{list}}
\newcounter{iii}
\newenvironment{iii}{\begin{list}{\rm \roman{iii}) }{\usecounter{iii} \leftmargin=0.0pt \labelsep=0.0pt \listparindent=0.0pt \labelwidth=0.0pt \parsep=\smallskipamount \itemsep=0.0pt \topsep=0.0pt \partopsep=\smallskipamount}}{\end{list}}

\newcommand{\bP}{\mathop{\text{\bf P}}\nolimits}
\newcommand{\Gal}{\mathop{\text{\rm Gal}}\nolimits}

\newcommand{\N}{\mathop{\text{\rm N}}\nolimits}

\newcommand{\bbC}{{\mathbbm C}}
\newcommand{\bbQ}{{\mathbbm Q}}
\newcommand{\bbZ}{{\mathbbm Z}}

\newcommand{\bra}{ }
\newcommand{\brr}{, }

\begin{document}

%
%
%
%
%
%
%
%
%

\title[On cubic surfaces with a rational line]
 {On cubic surfaces with a rational line}

\author{Andreas-Stephan Elsenhans}

\address{%
Mathematisches Institut\\
Universit\"at Bayreuth\\
Universit\"atsstra\ss e 30\\
D-95440 Bayreuth\\
Germany}

\email{Stephan.Elsenhans@uni-bayreuth.de}

\thanks{The first author was supported in part by the Deutsche Forschungsgemeinschaft (DFG) through a funded research~project.}
\author{J\"org Jahnel}
\address{D\'epartement Mathematik\\
Universit\"at Siegen\\
Walter-Flex-Str.~3\\
D-57068 Siegen\\
Germany}

\email{jahnel@mathematik.uni-siegen.de}

\subjclass{Primary 14J26; Secondary 14G25, 11G35}

\keywords{Cubic surface,
Degree~$4$
Del Pezzo surface, Explicit Galois descent}

\date{May 15, 2011}
\dedicatory{~}

\begin{abstract}
We report on our project to construct non-singular cubic surfaces
over~$\bbQ$
with a rational~line. Our~method is to start with
degree~$4$
Del Pezzo surfaces in diagonal~form. For~these, we develop an explicit version of Galois~descent.
\end{abstract}

\maketitle

\section{Introduction}

\begin{ttt}
The~configuration of the 27~lines upon a smooth cubic surface is highly~symmetric. The~group of all permutations respecting the canonical class as well as the intersection pairing is isomorphic to the Weyl
group~$W(E_6)$
of
order~$51\,840$.

When~$S$
is a cubic surface
over~$\bbQ$,
the absolute Galois
group~$\Gal(\overline\bbQ/\bbQ)$
operates on the 27~lines. This~yields a
subgroup~$G \subseteq W(E_6)$.
It~is an open problem whether each of the
$350$
conjugacy classes of subgroups
of~$W(E_6)$
may be realized by a cubic surface
over~$\bbQ$.

Exactly~$172$
of the
$350$
conjugacy classes fix a~line. We~constructed examples of cubic surfaces
over~$\bbQ$
realizing each of these~subgroups. The~goal of this note is to report on our~investigations.
\end{ttt}

\begin{rem}
The~analogous question for Del Pezzo surfaces of
degree~$4$
is somewhat easier as it leads to subgroups
of~$W(D_5)$.
B.\,\`E.~Kunyavskij, A.\,N.~Skorobogatov, and M.\,A.~Tsfasman~\cite{KST} showed that every subgroup
of~$W(D_5)$
may be realized by a surface defined
over~$\bbQ$.
\end{rem}

\section{Constructions}

\begin{ttt}
Cubic~surfaces with a rational line are closely related to Del Pezzo surfaces of
degree~$4$.
Indeed,~blowing down the line leads to a
degree~$4$
Del Pezzo surface having a rational~point. On~the other hand, blowing up a rational point on a
degree~$4$
Del Pezzo surface yields a cubic surface with a rational~line. These~two constructions may easily be made~explicit.
\end{ttt}

\begin{conso}[{\rm Cubic surfaces versus Del Pezzo surfaces of degree~$4$}{}]
Let~a base
field~$K$
be fixed once and for~all.

\begin{iii}
\item
For~two linear forms
$l_0, l_1$,
suppose that the line
$l_0 = l_1 = 0$
is contained in the cubic
surface~$S$
given by
$F(x_0, \ldots, x_3) = 0$.
Then,~$F$
may be written as
$F = l_0 q_0 + l_1 q_1$
for quadratic forms
$q_0$
and~$q_1$.
The~corresponding 
degree~$4$
Del Pezzo
surface~$V$
is given by
$q_0 + l_1 x_4 = q_1 - l_0 x_4 = 0$.
\item
On~the other hand, let a Del Pezzo
surface~$V$
of
degree~$4$
be given by
$Q_0 (x_0, \ldots, x_4) = Q_1 (x_0, \ldots, x_4) = 0$.
If~$(0:0:0:0:1) \in V$
then
$Q_0$
and~$Q_1$
may be written as
$Q_0 = q_0 + l_0 x_4$
and
$Q_1 = q_1 + l_1 x_4$
for
$q_0, q_1$
quadratic forms and
$l_0, l_1$
linear forms in
$x_0, \ldots, x_3$,~only.
The~corresponding cubic
surface~$S$
is given by
$q_0 l_1 - q_1 l_0 = 0$.
\end{iii}
\end{conso}

\begin{rems}
\begin{abc}
\item
These~two constructions are inverse to each~other.
\item
One~may start construction~ii) as well with arbitrary generators of the pencil spanned by 
$Q_0$
and~$Q_1$.
\end{abc}
\end{rems}

\begin{fac}
Let\/~$A$
be a symmetric matrix  representing the quadratic
form\/
$q_0 |_{l_0 = 0}$.
If~the eigenvalues
of\/~$A$
are\/~$z_1, z_2, z_3$
then there is a symmetric matrix
representing\/~$Q_0$
with eigenvalues\/~$(-1), 1, z_1, z_2, z_3$.
\end{fac}

\begin{cor}
\begin{iii}
\item
In~particular,
$Q_0$~is
of rank\/
$<\!5$
if and only if\/
$q_0 |_{l_0=0}$
is of
rank\/~$<\!3$.
Hence,~the five degenerate quadratic forms in the pencil\/
$[Q_0, Q_1]$
are in one-to-one correspondence with the five tritangent planes through the line~considered.
\item
If~the eigenvalues of a symmetric matrix representing\/~$Q_0$
are\/
$0, Z_1, \ldots, Z_4$
then\/
``$l_0 = 0$\!''
is a tritangent plane
on\/~$S$.
The~conic, defined
by\/~$S$
on this plane, splits into two lines over the
field\/~$K(\sqrt{\mathstrut Z_1Z_2Z_3Z_4})$.
\end{iii}
\end{cor}

\begin{ex}
Consider~the case that
$Q_0 := a_0 x_0^2 + \ldots + a_4 x_4^2$
and
$Q_1 := b_0 x_0^2 + \ldots + b_4 x_4^2$
are diagonal forms over the
field~$K$.
Then,~the five tritangent planes correspond to the points
$((-b_i) : a_i) \in \bP^1$
as
$(-b_i Q_0 + a_i Q_1)$
is~degenerate.
The~conics split over the~fields
$$\textstyle K\big(\sqrt{\mathstrut\smash{\prod\limits_{j \neq i} (-b_i a_j + a_i b_j)}}\big)$$
for~$i = 0, \ldots, 4$.
Observe~that the product of the five radicands is a perfect~square.

On~the corresponding cubic surface, all 27~lines are defined~over
$$\textstyle L = K\big(\sqrt{\mathstrut\smash{\prod\limits_{j \neq 0} (-b_0 a_j + a_0 b_j)}}, \ldots, \sqrt{\mathstrut\smash{\prod\limits_{j \neq 4} (-b_4 a_j + a_4 b_j)}}\big) \, .$$
Indeed,~the subgroup
of~$W(E_6)$
stabilizing a line is clearly of
order~$51840/27 = 1920$.
It~is actually the semi-direct product
$T \rtimes S_5$,
where
$T \subset (\bbZ/2\bbZ)^5$
is the subgroup of
order~$16$
formed by the elements having an even number of components equal
to~$1$.
As~$\Gal(\overline\bbQ/L)$
stabilizes not only the five tritangent planes but also the lines on them, it must act through the trivial subgroup
of~$T \rtimes S_5$.
\end{ex}

\begin{cono}[{}{\rm Explicit Galois descent}]
\label{expl}
Let~$A$
be a commutative \'etale algebra of
degree~$5$
over~$\bbQ$
and
$\iota_0, \ldots, \iota_4 \colon A \to \bbC$
be the five~embeddings.

\begin{iii}
\item
For~general~$a, \!b \in \!A$,
the equations
$$\iota_0(a) x_0^2 + \cdots + \iota_4(a) x_4^2 = \iota_0(b) x_0^2 + \cdots + \iota_4(b) x_4^2 = 0$$
define a Del Pezzo
surface~$V$
of
degree~$4$
over~$\overline\bbQ$.
\item
Let~$l$
be a linear form in five variables with coefficients
in~$A$.
Then,~by symmetry, the quadratic forms
$\iota_0(a) (l^{\iota_0})^2 + \cdots + \iota_4(a) (l^{\iota_4})^2$
and
$\iota_0(b) (l^{\iota_0})^2 + \cdots + \iota_4(b) (l^{\iota_4})^2$
have rational~coefficients.
If~$l^{\iota_0}, \ldots, l^{\iota_4}$ 
are linearly independent then we have a Del Pezzo
surface~$V_0$
of
degree~$4$
over~$\bbQ$
such that its base change
to~$\overline\bbQ$
is isomorphic
to~$V$.
\end{iii}
\end{cono}

\begin{rems}
\begin{abc}
\item
This~construction is analogous to~\cite[Theorem~6.1]{EJ1}.
\item
The~five tritangent planes
on~$V_0$
correspond to the points
$((-\iota_i(b)) : \iota_i(a)) \in \bP^1$.
Hence,~the Galois operation on them is the same as that on the
embeddings~$\iota_i$.
\item
When~$a \neq 0$,
the conic on the tritangent plane corresponding
to~$((-\iota_i(b)) : \iota_i(a))$
splits into two lines over the~field
\begin{equation}
\label{kp}
\textstyle \bbQ\big(\iota_i(-b/a), \sqrt{\mathstrut\smash{\prod\limits_{j \neq i} (-\iota_i(b) \iota_j(a) + \iota_i(a) \iota_j(b))}}\big) \, .
\end{equation}
The~radicand may be rewritten as
$\N(a) \,\iota_i(a^3 \,\delta_{A/\bbQ}(-b/a))$,
where
$\delta_{A/\bbQ}$
denotes the different of an element
of~$A$.
\end{abc}
\end{rems}

\begin{ttt}
\label{diff}
Thus,~given a subgroup
$G \subseteq T \!\rtimes\! S_5$,
there is the following strategy to construct a cubic
surface~$S$
over~$\bbQ$
such that
$\Gal(\overline\bbQ/\bbQ)$
operates
via~$G$
on the 27~lines.\medskip

\noindent
{\bf Strategy.}
i)
Find~a number
field~$K$, normal
over~$\bbQ$,
such that
$\Gal(K/\bbQ) \cong G$.
Identify~the normal subextension
$K^\prime \subseteq K$
such that
$\Gal(K^\prime/\bbQ)$
is the
image~$G^\prime$
of~$G$
in~$S_5$~\cite{MM}.

\begin{iii}
\addtocounter{iii}{1}
\item
Find~five elements
$r_0, \ldots, r_4 \in K^\prime$
with the properties~below.

$r_0, \ldots, r_4$
are permuted
by~$G^\prime$
exactly via the embedding
$G^\prime \subseteq S_5$.
Further,~the square roots
$\pm\sqrt{\mathstrut r_0}, \ldots, \pm\sqrt{\mathstrut r_4}$
are elements
of~$K$
and acted upon
by~$G$
according to the embedding
$G \subseteq T \!\rtimes\! S_5$.

Put~$p(T) := (T - r_0) \cdot\ldots\cdot (T - r_4)$
and~$A := \bbQ[T]/(p)$.
This~is a commutative \'etale algebra of
degree~$5$
over~$\bbQ$
with a distinguished element
\mbox{$r := (T \!\!\mod (p))$}.
\item
Choose~$x \in A$
and put
$d := \delta_{A/\bbQ}(x)$.
Set~$a := dr$
and~$b := -xa$.
\item
Execute~Construction~\ref{expl}
for~$a,b \in A$.
On~the Del Pezzo
surface~$V_0$
found, search for a
\mbox{$\bbQ$-rational}~point.
If~none is found then go back to step~iii). Otherwise,~determine the cubic
surface~$S$.
\end{iii}
\end{ttt}

\begin{rems}
\begin{abc}
\item
The~properties required in~ii) imply
$\sqrt{\mathstrut r_0} \cdot\ldots\cdot \sqrt{\mathstrut r_4} \in \bbQ$.
I.e.,~$\N(r)$
is a perfect~square.
\item
The~construction yields
$\N(a)a^3\,\delta_{A/\bbQ}(-b/a) = \N(a)(d^2r)^2r$.
As~the product of the five radicands
in~(\ref{kp})
is a square, the norm
of~$a$
is a perfect square~automatically.
\end{abc}
\end{rems}

\section{Examples}

\begin{ttt}
There~are 172 conjugacy classes of subgroups
of~$W(E_6)$
that fix a~line. We~constructed examples for each such~group.

Actually,~$81$
of the
$172$~classes
also stabilize a double-six and
$49$
of the
$172$~classes
stabilize a pair of Steiner~trihedra.
$34$~classes
do~both.
Thus,~examples for
$96$
of the
$172$~conjugacy
classes had been constructed before~\cite{EJ1,EJ2}.
The~remaining
$76$~classes
were of~interest.

After~naive trials and an extensive search through surfaces with small coefficients, only six of the
$76$
classes remained~open. For~these, we applied Strategy~\ref{diff}.
\end{ttt}

\begin{rem}
In~Strategy~\ref{diff}, we regularly run into reiteration, because there were no
\mbox{$\bbQ$-rational}
points on the Del Pezzo surfaces of
degree~$4$.
\end{rem}

\begin{ttt}
The~list containing our examples of cubic surfaces is available on the second author's web page at {\tt http://www.uni-math.gwdg.de/jahnel/\discretionary{}{}{}Arbeiten/Kub\_Fl/list\_rat\_ger.txt}. The~numbering of the subgroups is that created by {\tt GAP}, version~4.4.12.
\end{ttt}

\begin{ex}
As~a conclusion, let us show how Strategy~\ref{diff} works on a particular~example. We~consider the subgroup of number~107.

Abstractly,~this is a
group~$G$
of
order~$16$.
Its~center is isomorphic to the Klein four-group. The~operation on the 27~lines causes orbits of lengths
$1$,
$2$,
$4$,
$4$,
and~$16$.
On~the two orbits of size four,
$G$~acts
via two different quotients, both isomorphic to the dihedral
group~$D_4$
of order~eight. The~operation on the five tritangent planes through the rational line is via a
quotient~$G^\prime$
of order~four. The~orbits are of sizes
$1$,
$2$,
and~$2$.

\begin{iii}
\item
An~example of a field with Galois
group~$G$
is the
composite~$K := K_1K_2$
of
$\smash{K_1 := \bbQ\big(\!\sqrt{3 \pm \sqrt{3}}\big)}$
and~$\smash{K_2 := \bbQ\big(\!\sqrt{-9 \pm \sqrt{6}}\big)}$.
Then,~the subfield corresponding
to~$G^\prime$
is~$K^\prime = \bbQ(\sqrt{3}, \sqrt{6}) = \bbQ(\sqrt{2}, \sqrt{3})$.
Observe~that both fields
$K_1$
and~$K_2$
contain~$K^\prime$.
Further,~both are extensions
of~$\bbQ$
of
type~$D_4$.
\item
Thus,~we chose
$r_0, \ldots, r_4$
to be
$2, 3 \pm \sqrt{3}$,
and~$-9 \pm \sqrt{6}$.
This~yields
$$p(T) = (T - 2)[(T-3)^2 - 3][(T+9)^2 - 6] \, .$$
\item
We~worked with
$x := r = (T \!\!\mod (p))$.
\item
The biggest coefficient of the resulting del Pezzo
surface~$V_0$
is
$524\,391\,211\,895\,464$.
An~isomorphic surface is given by the~equations
\begin{eqnarray*}
 &   & 4 x_0^2 + 10 x_0 x_1 + 20 x_0 x_2 - 112 x_0 x_3 - 134 x_0 x_4 + 7 x_1^2 - 26 x_1 x_2 - 134 x_1 x_3 \\
 &   & {} - 148 x_1 x_4 - 2 x_2^2 + 140 x_2 x_3 - 2 x_2 x_4 + 10 x_3^2 - 38 x_3 x_4 - 323 x_4^2 \\
 & = & 47 x_0^2 - 18 x_0 x_1 + 10 x_0 x_2 - 188 x_0 x_3 - 178 x_0 x_4 \!+\! 63 x_1^2 \!-\! 22 x_1 x_2 \!+\! 376 x_1 x_3 \\
 &   & {} - 86 x_1 x_4 + 71 x_2^2 - 580 x_2 x_3 + 146 x_2 x_4 - 364 x_3^2 - 296 x_3 x_4 - 21 x_4^2 = 0 \, .
\end{eqnarray*}
Here,~a point search in {\tt magma} with an initial height limit of 100 shows 14 rational~points. Blowing~up
$(8 : -13 : 4 : 2 : -3)$
leads to a cubic surface with coefficients up
to~$3\,838\,320$.
Reembedding~gives us the final result, the cubic
surface~$V$
with the~equation
\begin{eqnarray*}
& & 2 x^2 y + 6 x^2 z - 4 x y^2 + 6 x y z + 4 x y w - 10 x z^2 - 4 x z w -
7 x w^2 + 2 y^3 \\
& & {} - 9 y^2 z - 4 y^2 w + 4 y z^2 - 26 y z w + 6 y w^2 +
z^3 + 10 z^2 w - 7 z w^2 - 5 w^3 = 0 \, .
\end{eqnarray*}
\end{iii}
\end{ex}

\begin{rem}
The~rational line
on~$V$
connects
$(5 : 0 : 0 : -7)$
with
$(0 : 5 : 10 : 2)$.
\end{rem}

\begin{rem}
There~are actually a few more particularities characterizing the subgroup of number~107.

\begin{abc}
\item
First~of all, the two
$D_4$~extensions
$K_1$
und~$K_2$
become cyclic over the same quadratic
field~$\bbQ(\sqrt{2})$.
\item
On~the other hand,
over~$\bbQ(\sqrt{3})$
and~$\bbQ(\sqrt{6})$,
they are of Kleinian~type. However,~there is yet another~oddity.
While~$\Gal(K_1/\bbQ(\sqrt{3}))$
operates on the corresponding four lines via two disjoint two-cycles,
$\Gal(K_2/\bbQ(\sqrt{3}))$
acts on its orbit by double-trans\-po\-si\-tions.
Over~$\bbQ(\sqrt{3})$
instead
of~$\bbQ(\sqrt{6})$,
the situation is vice~versa.
\end{abc}
To~realize such a behaviour, it was essential to choose
$r_1$
in
$\bbQ(\sqrt{3})$
fulfilling
$\N(r_1) \in 6(\bbQ^*)^2$
and
$r_3$
in
$\bbQ(\sqrt{6})$
such that
$\N(r_3) \in 3(\bbQ^*)^2$.
\end{rem}


\end{document}